\newcommand{\al}{\alpha}

\newcommand{\veps}{\varepsilon}          

\newcommand{\cal}{\mathcal}

\newcommand{\calv}{{\cal V}}

\newcommand{\Fix}{{\rm Fix}}

\newcommand{\incl}{\subseteq}

          \newcommand{\limpl}{\Longrightarrow}
      
\newcommand{\oo}{\infty}

\newcommand{\bk}{\bigskip}               \newcommand{\sk}{\smallskip}
         \newcommand{\n}{\noindent}

                \def\R+oo{R_+\cup\{\oo\}}

\def\dtends   {\stackrel {\it d}{\longrightarrow}}
\def\Dtends   {\stackrel {\it D}{\longrightarrow}}
\def\etends   {\stackrel {\it e}{\longrightarrow}}

\def\(V)tends {\stackrel {(\calv)}{\longrightarrow}}

\newcommand{\barr}{\begin{array}}          \newcommand{\earr}{\end{array}}
\newcommand{\beq}{\begin{equation}}        \newcommand{\eeq}{\end{equation}}
\newcommand{\bit}{\begin{itemize}}         \newcommand{\eit}{\end{itemize}}
\newcommand{\blemma}{\begin{lemma}}        \newcommand{\elemma}{\end{lemma}}
\newcommand{\bprop}{\begin{proposition}}   \newcommand{\eprop}{\end{proposition}}
\newcommand{\bproof}{\begin{proof}}        \newcommand{\eproof}{\end{proof}}
\newcommand{\btab}{\begin{tabular}}        \newcommand{\etab}{\end{tabular}}
\newcommand{\btheorem}{\begin{theorem}}    \newcommand{\etheorem}{\end{theorem}}

\documentclass[reqno]{amsart}

\newtheorem{theorem}{\bf Theorem}

\newtheorem{lemma}{\bf Lemma}
\newtheorem{proposition}{\bf Proposition}

\begin{document}

\title
[Nieto-Lopez Theorems in Ordered Metric Spaces]
{NIETO-LOPEZ THEOREMS \\
IN ORDERED METRIC SPACES}

\author{Mihai Turinici}
\address{
"A. Myller" Mathematical Seminar;
"A. I. Cuza" University;
11, Copou Boulevard; 
700506 Ia\c{s}i, Romania
}
\email{mturi@uaic.ro}




\begin{abstract}
The comparison type version of the fixed point result 
in ordered metric spaces established by Nieto and Rodriguez-Lopez
[Acta Math. Sinica (English Series), 23 (2007), 2205-2212]
is nothing but a particular case of the classical
Banach's contraction principle
[Fund. Math., 3 (1922), 133-181].
\end{abstract}

\maketitle

\section{Introduction}
\setcounter{equation}{0}

Let $(X,d) $ be a metric space; and $T:X\to X$ be a selfmap of $X$.
We say that $x\in X$ is a {\it Picard point} (modulo $(d,T)$) if
{\bf i)} $(T^nx; n\ge 0)$ (=the {\it orbit} of $x$) is $d$-convergent, 
{\bf ii)} $z:=\lim_nT^nx$ is in $\Fix(T)$ (i.e.: $z=Tz$).
If this happens for each $x\in X$ and
{\bf iii)} $\Fix(T)$ is a singleton, 
then $T$ is referred to as a {\it Picard operator} (modulo $d$); 
see Rus \cite[Ch 2, Sect 2.2]{rus-2001}.

For example, such a property holds whenever 
$d$ is {\it complete} and
$T$ is {\it $d$-contractive}; cf. (b04).
A structural extension of this fact
-- when an {\it order} $(\le)$ on $X$ is being added --
was obtained in 2007 by
Nieto and Rodriguez-Lopez \cite{nieto-rodriguez-lopez-2007}.
Denote
\bit
\item[(a01)]
($x,y\in X$):\ $x<> y$ iff either $x\le y$ or $y\le x$
($x$ and $y$ are comparable).
\eit
This relation is reflexive and symmetric; 
but not in general transitive.
Let us say that the sequence $(x_n; n\ge 0)$ in $X$ is
$<>$-{\it ascending} if $x_n<>x_{n+1}$, for all $n$;
i.e.: any two consecutive terms of it are comparable.
Call the selfmap $T$, {\it $(d,\le;\al)$-contractive} (where $\al> 0$), if
\bit
\item[(a02)]
$d(Tx,Ty)\le \al d(x,y)$,\ $\forall x,y\in X$, $x\le y$.
\eit
If this holds for some $\al\in ]0,1[$, 
we say that $T$ is {\it $(d,\le)$-contractive}.

\btheorem \label{t1}
Assume that $d$ is complete, $T$ is $(d,\le)$-contractive,
and
\bit
\item[(a03)]
$X(T,<>):=\{x\in X; x<> Tx\}$ is nonempty
\item[(a04)]
$T$ is monotone (increasing or decreasing)
\item[(a05)]
for each $x,y\in X$, $\{x,y\}$ has lower and upper bounds
\item[(a06)]
each $<>$-ascending sequence $(x_n; n\ge 0)$ with $x_n\dtends x$ \\
has a subsequence $(y_n:=x_{q(n)}; n\ge 0)$
with $y_n<> x$, $\forall n$.
\eit
Then, $T$ is a Picard operator (modulo $d$).
\etheorem

Note that, this conclusion is retainable as well
when (a06) is replaced with
\bit
\item[(a07)]
$T$ is $d$-continuous: $x_n\dtends x$ $\limpl$ $Tx_n\dtends Tx$;
\eit
this is just the 2004 main result in 
Ran and Reurings \cite{ran-reurings-2004}.
According to many authors, 
these two results are credited to be the first extension of the classical 1922
Banach's contraction mapping principle \cite{banach-1922}
to the realm of (partially) ordered metric spaces.
Unfortunately, the assertion is not true:
some early statements of this type 
have been obtained two decades ago by
Turinici \cite{turinici-1986},
in the context of ordered metrizable uniform spaces.
(We refer to Section 4 below for details).
\sk \sk

Now, the Nieto-Rodriguez-Lopez fixed point
result found some useful applications to 
differential and integral equations theory; cf.
O'Regan and Petru\c{s}el \cite{o-regan-petrusel-2008}.
So, it cannot be surprising that, soon after, many extensions of
Theorem \ref{t1} were provided;
for the most consistent contributions we refer to
Agarwal, El-Gebeily and O'Regan
\cite{agarwal-el-gebeily-o-regan-2008},
Gwozdz-Lukawska and Jachymski
\cite{gwozdz-lukawska-jachymski-2009},
or
Ciric et al \cite{ciric-mihet-saadati-2009}.
It is therefore natural to discuss
the position of Theorem \ref{t1}
within the classification scheme proposed by
Rhoades \cite{rhoades-1977};
see also 
Collaco and E Silva \cite{collaco-e-silva-1997}.
The conclusion to be derived reads (cf. Section 2):
the Nieto-Rodriguez-Lopez theorem is but a particular case of the 
Banach's contraction principle \cite{banach-1922}.
Further, in Section 3, a 
Suzuki type variant  \cite{suzuki-2008}
of Theorem \ref{t1} is considered.
Some other aspects will be delineated elsewhere.

\section{Main result}
\setcounter{equation}{0}

Let again $(X,\le,d)$ be an ordered metric space; 
and $T:X\to X$, a selfmap of $X$.
Given $x,y\in X$, any subset 
$\{z_1,...,z_k\}$ (for $k\ge 2$) in $X$ with
$z_1=x$, $z_k=y$, and [$z_i<> z_{i+1}$,  $i\in \{1,...,k-1\}$]
will be referred to as a {\it $<>$-chain} 
between $x$ and $y$; the class of all these will 
be denoted $C(x,y;<>)$.
Let $\sim$ stand for the relation over $X$ attached to $<>$ as 
\bit
\item[(b01)]
$x\sim y$ iff $C(x,y;<>)$ is nonempty. 
\eit
Clearly, $(\sim)$ is reflexive and symmetric; because so is 
$<>$. Moreover, $(\sim)$ is transitive; 
hence, it is an equivalence over $X$.

The following variant of Theorem \ref{t1} is our starting point.

\btheorem \label{t2}
Assume that $d$ is complete, 
$T$ is $(d,\le)$-contractive, condition (a06) holds, and
\bit
\item[(b02)]
$T$ is $<>$-increasing [$x<> y$ implies $Tx<> Ty$]
\item[(b03)]
$(\sim)=X\times X$
[$C(x,y; <>)$ is nonempty, for each $x,y\in X$].
\eit
Then, $T$ is a Picard operator (modulo $d$).
\etheorem

This result includes Theorem \ref{t1};
because (a04) $\limpl$ (b02), (a05) $\limpl$ (b03).
[For, given $x,y\in X$,  
there exist, by (a05), some $u,v\in X$ with
$u\le x\le v$, $u\le y\le v$. This yields $x<>u$, $u<> y$; 
wherefrom, $x\sim y$]. In addition, it tells us
that the regularity condition (a03) is not needed there.

Now, the remarkable fact to be underlined is that Theorem \ref{t2}
(hence the Nieto-Rodriguez-Lopez statement as well) 
is deductible from the 1922
Banach's contraction mapping principle \cite{banach-1922}.
Let $e(.,.)$ be another metric over $X$.
Call $T:X\to X$, 
{\it $(e;\al)$-contractive} (for some $\al> 0$) when 
\bit
\item[(b04)]
$e(Tx,Ty)\le \al e(x,y)$, $\forall x,y\in X$; 
\eit
if this holds for some $\al\in ]0,1[$, 
the resulting convention will read as: 
$T$ is {\it $e$-contractive}.
The announced Banach's result is:

\btheorem \label{t3}
Assume that $e$ is complete and $T$ is $e$-contractive.
Then, $T$ is a Picard operator (modulo $e$).
\etheorem

We are now in position to give the announced answer.

\bprop \label{p1}
We have
Theorem \ref{t3} $\limpl$ Theorem \ref{t2};
hence (by the above) 
the Banach fixed point principle implies the Nieto-Rodriguez-Lopez
result.
\eprop

\bproof
Let the conditions of Theorem \ref{t2} hold.
We introduce a mapping $e:X\times X\to R_+$ as:
for each $x,y\in X$,
\bit
\item[(b05)]
$e(x,y)=\inf[d(z_1,z_2)+...+d(z_{k-1},z_k)]$, \\
where $\{z_1,...,z_k\}$ is a $<>$-chain between $x$ and $y$.
\eit

{\bf I)}
Clearly, $e$ is 
reflexive [$e(x,x)=0$, $\forall x\in X$],
symmetric [$e(y,y)=e(y,x)$, $\forall x,y\in X$]
and 
triangular [$e(x,z)\le e(x,y)+e(y,z)$, $\forall x,y,z\in X$].
In addition, the triangular property of $d$ gives 
$d(x,y)\le d(z_1,z_2)+...+d(z_{k-1},z_k)$,
for any $<>$-chain $\{z_1,...,z_k\}$  
between $x$ and $y$. So, passing to infimum, yields 
\beq \label{201}
\mbox{
$d(x,y)\le e(x,y)$, $\forall x,y\in X$
[referred to as: $d$ is subordinated to $e$].
}
\eeq
Note that 
$e$ is sufficient in such a case [$e(x,y)=0$ $\limpl$ $x=y$];
hence, it is a (standard) metric on $X$.
Finally, by the very definition of $e$, one has 
\beq \label{202}
\mbox{
$d(x,y)\ge e(x,y)$ (hence $d(x,y)=e(x,y)$), whenever $x<>y$.
}
\eeq

{\bf II)}
We claim that $e$ is complete on $X$. 
Let $(x_n; n\ge 0)$ be an $e$-Cauchy sequence in $X$.
There exists a strictly ascending sequence of ranks
$(k(n); n\ge 0)$, in such a way that
[($\forall n$): $k(n)< m$ $\limpl$ $e(x_{k(n)},x_m)< 2^{-n}$].
Denoting $(y_n:=x_{k(n)}, n\ge 0)$, we therefore have
$e(y_n,y_{n+1})< 2^{-n}$, $\forall n$.
Moreover, by the imposed $e$-Cauchy property, 
$(x_n)$ is  $e$-convergent iff so is $(y_n)$. 
To establish this last property, one may proceed as follows.
As $e(y_0,y_1)< 2^{-0}$, there exists 
(for the starting  rank $p(0)=0$)
a $<>$-chain $\{z_{p(0)},...,z_{p(1)}\}$ between $y_0$ and $y_1$
(hence $p(1)-p(0)\ge 1$, $z_{p(0)}=y_0$, $z_{p(1)}=y_1$), 
such that 
$$
d(z_{p(0)},z_{p(0)+1})+...+d(z_{p(1)-1},z_{p(1)})< 2^{-0}.
$$
Further, as $e(y_1,y_2)< 2^{-1}$, there exists a 
$<>$-chain $\{z_{p(1)},...,z_{p(2)}\}$ between $y_1$ and $y_2$
(hence $p(2)-p(1)\ge 1$, $z_{p(1)}=y_1$, $z_{p(2)}=y_2$), 
such that 
$$
d(z_{p(1)},z_{p(1)+1})+...+d(z_{p(2)-1},z_{p(2)})< 2^{-1};
$$
and so on. The procedure may continue indefinitely;
it gives us a $<>$-ascending sequence $(z_n; n\ge 0)$ in $X$ with 
(cf. (\ref{202}))
\beq \label{203}
\sum_n e(z_n,z_{n+1})=\sum_n d(z_n,z_{n+1})< \sum_n 2^{-n}< \oo.
\eeq
In particular, $(z_n; n\ge 0)$ is $d$-Cauchy; 
wherefrom (as $d$ is complete),
$z_n\dtends z$ as $n\to \oo$, for some $z\in X$.
Combining with (a06), there must be a subsequence 
$(t_n:=z_{q(n)}; n\ge 0)$ of $(z_n; n\ge 0)$ with
$t_n<> z$, $\forall n$. This firstly gives 
(by the previous relation), $t_n\dtends z$ as $n\to \oo$.
Secondly (again via (\ref{202})),
$e(t_n,z)=d(t_n,z)$, $\forall n$; so that 
[combining with the above $d$-convergence property], 
$t_n\etends z$ as $n\to \oo$; 
On the other hand, (\ref{203}) also tells us that
$(z_n; n\ge 0)$ is $e$-Cauchy. 
Adding the $e$-convergence property of $(t_n; n\ge 0)$ gives 
$z_n\etends z$ as $n\to \oo$;
wherefrom (as $z_{p(n)}=y_n, n\ge 0$), 
$y_n\etends z$ as $n\to \oo$; and our claim follows.

{\bf III)}
Let $\al\in ]0,1[$ be the number appearing in 
the $(d,\le)$-contractive property of $T$. 
Given $x,y\in X$, let $\{z_1,...,z_k\}$ be 
a $<>$-chain connecting them (existing via (b03)).
From (b02), $\{Tz_1,...,Tz_k\}$ 
is a $<>$-chain between $Tx$ and $Ty$. 
So, combining with the contractive condition,
$$
e(Tx,Ty)\le \sum_{i=1}^{k-1} d(Tz_i,Tz_{i+1})\le 
\al \sum_{i=1}^{k-1}d(z_i,z_{i+1}),
$$
for all such $<>$-chains; wherefrom, passing to infimum, 
$e(Tx,Ty)\le \al e(x,y)$; i.e., (b04) holds.
Summing up, Theorem \ref{t3} applies to these data; and we are done.
\eproof

\section{Suzuki type extensions}
\setcounter{equation}{0}

{\bf (A)}
Let $(X,d)$ be a metric space.
Define a mapping $F:[0,\oo[\to ]0,1]$ as
\bit
\item[(c01)]
$F(t)=1$, if $0\le t\le (\sqrt{5}-1)/2$ \\
$F(t)=(1-t)t^{-2}$, if $(\sqrt{5}-1)/2\le t\le 2^{-1/2}$ \\
$F(t)=(1+t)^{-1}$, if $2^{-1/2}\le t< \oo$.
\eit
Note that $F$ is continuous and decreasing over its existence domain.
Given the selfmap $T:X\to X$, call it {\it conditional $(F,\al)$-contractive}
(where $\al\ge 0$) provided
\bit
\item[(c02)]
[$x,y\in X$, $F(\al)d(x,Tx)\le d(x,y)$] $\limpl$ $d(Tx,Ty)\le \al d(x,y)$;
\eit
if this holds for some $\al\in [0,1[$, then $T$ is called 
{\it conditional $F$-contractive}.
The following 2008 result in Suzuki \cite{suzuki-2008} 
is our starting point.

\btheorem \label{t4}
Suppose that $d$ is complete and $T$ is conditional 
$F$-contractive. Then, $T$ is a Picard operator.
\etheorem

It is worth remarking that the family of conditions (c02) gives a 
characterization of completeness;
see the quoted paper for details.
A related statement involving the 
Kannan type conditions \cite{kannan-1968}
is to be found in 
Kikkawa and Suzuki \cite{kikkawa-suzuki-2008}.
For various extensions of such results we refer to
Altun and Erduran \cite{altun-erduran-2011};
see also
Popescu  \cite{popescu-2009}.
Note that, in all these statements, 
the premise of the conditional contractive property (c02)
is "asymmetric" with respect to the couple $(x,y)$;
so, it is natural to ask whether
a supplementary condition may be added there,
with a "dual" information about the variable $y$.
It is our aim in the following to show that a positive answer to 
this is possible, in a (quasi-) ordered realm.
\sk

{\bf (B)}
Let $(X,d)$ be a metric space.
Take a quasi-order $(\le)$ 
(i.e.: reflexive and transitive relation) over $X$; 
as well as a selfmap $T:X\to X$.
Term the sequence $(z_n)$, 
{\bf i)} {\it ascending} iff $z_i\le z_j$ when $i\le j$,
{\bf ii)} {\it orbital}, when $z_n=T^nx$, $n\ge 0$, for some $x\in X$.
Call the ambient metric $d$, {\it ascending-orbital complete}
(in short: ao-complete) when each ascending orbital $d$-Cauchy sequence
converges.
Further, let us say that $(\le)$ is {\it ascending-orbital-self-closed}
(in short: ao-self-closed) provided:
[$(z_n)$=ascending orbital, $z_n\to z$] imply $z_n\le z$, for all $n$.
\sk

Denote
$G(t)=1/(1+t)$, $t> 0$;
this function is continuous, decreasing 
and maps $]0,\oo[$ onto  $]0,1[$.
Further, let the relation $<>$ over $X$ be introduced as in (a01), 
but, in our quasi-order realm;
as before, it  is reflexive and symmetric.
Call the selfmap $T$, {\it weakly conditional $(G,<>;\al)$-contractive}
(where $\al > 0$) provided
\bit
\item[(c03)]
$x<>y$, $G(\al)\max\{d(x,Tx),d(y,Tx)\}\le d(x,y)$ 
$\limpl$ $d(Tx,Ty)\le \al d(x,y)$.
\eit
If this holds for at least one $\al\in ]0,1[$, the underlying map $T$ is called
{\it weakly conditional $(G,<>)$-contractive}.

Having these precise, assume in the sequel that
\bit
\item[(co4)]
$X(T,\le):=\{x\in X; x\le Tx\}$ is nonempty
\item[(c05)]
$T$ is increasing ($x\le y$ $\limpl$ $Tx\le Ty$).
\eit
We say that $x\in X(T,\le)$ is a {\it Picard point} (modulo $(\le,T)$) if
{\bf iii)} $(T^nx; n\ge 0)$ is $d$-convergent, 
{\bf iv)} $z:=\lim(T^nx)$ is in $\Fix(T)$ and $T^nx\le z$, $\forall n$.
If this happens for each $x\in X(T,\le)$ and
{\bf v)} $\Fix(T)$ is $(\le)$-{\it singleton} [$z,w\in \Fix(T)$, $z\le w$ $\limpl$ $z=w$],
then $T$ is called a {\it Picard operator} (modulo $(\le)$). 
Note that, in this case, each $x^*\in \Fix(T)$ fulfills
\beq \label{301}
\forall u\in X(T,\le):\ x^*\le u \limpl u\le x^*;
\eeq
i.e.: $x^*$ is $(\le)$-maximal in $X(T,\le)$. 
In fact, assume that $x^*\le u\in X(T,\le)$. 
By iii) and iv), $(T^nu; n\ge 0)$
$d$-converges to some $u^*\in \Fix(T)$ with 
$T^nu\le u^*$, $\forall n$; hence, $x^*\le u\le u^*$.
Combining with v) gives  $x^*=u^*$; 
wherefrom $u\le x^*$.

\btheorem \label{t5}
Assume that
$T$ is weakly conditional $(G,<>)$-contractive,
$d$ is ao-complete, and $(\le)$ is ao-self-closed.
Then $T$ is a  Picard operator (modulo $(\le)$). 
\etheorem

\bproof
Let $\al\in ]0,1[$ be the number appearing in the weak 
conditional contractive property of $T$.
There are several steps to be passed.

{\bf I)}
Let $x^*,u^*\in \Fix(T)$ be such that $x^*\le u^*$.
We have
$d(x^*,Tx^*)=0$, $d(u^*,Tx^*)=d(u^*,x^*)$;
hence $G(\al)\max\{d(x^*,Tx^*),d(u^*,Tx^*)\}\le d(x^*,u^*)$.
This, by the contractive condition, yields
$d(x^*,u^*)\le \al d(x^*,u^*)$; wherefrom 
(as $d$=metric) $x^*=u^*$;
and so, $\Fix(T)$ is $(\le)$-singleton.

{\bf II)}
It remains to show that each $x=x_0\in X(T,\le)$ is a Picard point
(modulo $(\le,T)$).
Put $x_n=T^nx$, $n\ge 0$; 
clearly, $(x_n; n\ge 0)$ is an ascending orbital sequence in $X$.
For each $n\ge 0$,
$$
x_n\le x_{n+1}, 
G(\al)\max\{d(x_n,Tx_n), d(x_{n+1},Tx_n)\}\le d(x_n,x_{n+1});
$$
so, by the imposed contractive condition,
$d(x_{n+1},x_{n+2})\le \al d(x_n,x_{n+1})$, $\forall n$;
wherefrom, $(x_n; n\ge 0)$ is $d$-Cauchy.
As $d$ is ao-complete (and $(\le)$ is ao-self-closed)
\beq \label{302}
\mbox{
$x_n\to z$ (hence $x_n\le z$, $\forall n$), for some $z\in X$.
}
\eeq

{\bf III)}
Suppose that our sequence is such that
\bit
\item[(c06)]
for each $n$, there exists $m> n$ with $x_m=z$.
\eit
It follows that a subsequence $(x_{p(n)}; n\ge 0)$ of $(x_n)$ 
exists with $x_{p(n)}=z$, for all $n$. 
This, along with [$x_{p(n)+1}=Tz$, $\forall n$]  gives 
(via (\ref{302}) and $d$=metric), $z=Tz$.

{\bf IV)}
Assume in the following that the opposite situation holds:
\bit
\item[(c07)]
there exists $h$ such that:\ \ $x_n\ne z$, for all $n\ge h$.
\eit
Fix $k\ge h$; and put $x_k=u$; clearly,
$u\le x_n\le z$, for all $n\ge k$.
As 
$d(x_n,Tx_n)\to 0$, $d(x_n,u)\to d(z,u)> 0$, 
there must be some rank $p\ge k$ such that
$d(x_n,Tx_n)\le \al d(x_n,u)$, $\forall n\ge p$.
On the other hand (for the same ranks)
$d(u,Tx_n)\le d(u,x_n)+d(x_n,Tx_n)\le (1+\al) d(u,x_n)$;
hence, summing up,
$$
x_n\ge u,
G(\al) \max\{d(x_n,Tx_n),d(u,Tx_n)\}\le d(x_n,u),\ \
\forall n\ge p.
$$
These, by the contractive condition, give 
$d(Tx_n,Tu)\le \al d(x_n,u)$, $\forall n\ge p$;
so that (passing to limit as $n\to \oo$),
$d(z,Tu)\le \al d(u,z)$.
By the triangle inequality
$d(u,Tu)\le d(u,z)+d(z,Tu)\le (1+\al) d(u,z)$;
wherefrom (putting these together)
$$
u\le z,
G(\al) \max\{d(u,Tu),d(z,Tu)\}\le d(u,z);
$$
so that (by the same contractive condition),
$d(Tu,Tz)\le \al d(u,z)$.
Taking into account the adopted notation, we have
$d(Tx_k,Tz)\le \al d(x_k,z)$, $\forall k\ge h$.
So, passing to limit as $k\to \oo$, one derives 
$z=Tz$; and conclusion follows.
\eproof

In particular, when $(\le)$ is the trivial quasi-order of $X$,
the obtained result extends, in a partial way, Theorem \ref{t4}.
An open question is of whether or not
a full extension may be reached; 
we conjecture that a positive answer is not in general valid.
Note that, by the same technique, it is possible
to get a quasi-ordered  version of the main result in
Singh, Pathak and Mishra \cite{singh-pathak-mishra-2010};
we shall develop such questions elsewhere.

\section{Old approach (1986)}
\setcounter{equation}{0}

In the following, a summary of the 1986 general results in 
Turinici \cite[Sect 2-3]{turinici-1986} 
is being sketched, for completeness reasons.
\sk

{\bf (A)}
Let $X$ be a nonempty set; and $\le$ be an {\it ordering} 
(i.e., a reflexive, antisymmetric, and transitive relation) on $X$.
For any $x\in X$ denote $(\le,x]=\{y\in X; y\le x\}$ and $[x,\le)=\{y\in X; x\le y\}$; 
also, given any couple $x,y\in X$, $x\le y$, put $[x,y]=(\le,y] \cap [x,\le)$ 
and call it the (order) {\it interval} between $x$ and $y$.
A sequence $(x_n;n\in N)$ in $X$ will be said to be {\it increasing} 
when $x_i\le x_j$ for $i\le j$, and {\it bounded from above} in case $x_n\le  y$, $n\in N$,
for some $y\in X$.
Furthermore, let $D=(d_i;i\in N)$  be a denumerable sufficient family 
of {\it semi-metrics} on $X$; in this case, the triplet $(X,D,\le)$
will be termed an {\it ordered metrizable uniform space}. 
We shall say that the sequence $(x_n; n\in N)$ in $X$, $D$-{\it converges} to $x\in X$ 
(and we write $x_n\Dtends x$) when $d_i(x_n,x)\to 0$ as $n\to \oo$, for each $i\in N$;
if such elements exist, $(x_n; n\in N)$ will be called $D$-{\it convergent}. 
Further, let us say that $(x_n; n\in N)$ in $X$ is $D$-{\it Cauchy}  provided
$d_i(x_n,x_m)\to 0$  as $n,m\to \oo$, for each $i\in N$.
Clearly, any $D$-convergent sequence is necessarily $D$-Cauchy;
in this context, $(X,D)$ will be said to be {\it order complete} when 
each increasing $D$-Cauchy sequence converges.
A subset $Y$ of $X$ will be termed {\it order closed} when 
the limit of any $D$-convergent increasing sequence in $Y$ belongs to $Y$; 
also, the ambient ordering  on $X$ 
will be called {\it self-closed} (resp., {\it anti self-closed}) 
in case $[x,\le)$(resp., $(\le,x]$) 
is order closed for any $x$ in $X$; 
and {\it interval-closed}, when it is both self-closed  and anti self-closed 
(or, equivalently: when each interval of $X$ is order closed).

In what follows, we shall say that $(y_n; n \in N)$
is a {\it subsequence} of $(x_n; n \in N)$
when a strictly increasing function $k$ from $N$ to itself 
may be found with $x_{k(n)}=y_n$, $n\in N$. 
Under such a convention, let us call 
the sequence $(x_n; n\in N)$ in $X$,
{\it relatively compact} when any subsequence $(y_n; n\in N)$ of it 
contains a convergent subsequence. 
The following result,  closely related to that of 
Ward \cite{ward-1954} 
(see also 
Krasnoselskii \cite[Ch I, Sect 5]{krasnoselskii-1962})
will be useful for us.

\blemma \label{le1}
Let the ordered metrizable uniform space $(X,D;\le)$ 
be such that $\le$ is interval closed. 
Then, the increasing sequence $(x_n; n\in N)$ in $X$ is relatively compact, 
if and only if it converges to some element $x$ of $X$.
\elemma

\bproof
Let $(u_n; n\in N)$ and $(v_n; n\in N)$ be a couple of convergent subsequences 
of $(x_n; n\in N)$. If $u_n \Dtends u$ and $v_n \Dtends v$ then, by the 
interval-closeness property, we get $u\le v\le u$; that is, $u=v$. 
In other words, all convergent subsequences of $(x_n; n \in N)$ have the same limit, $x$.
We claim that $x_n \Dtends x$. Indeed, suppose that this assertion were false; then, a couple
$i\in N$, $\veps> 0$ may be chosen so that, for each $n\in N$ there exists $m> n$ with
$d_i(x_m,x)\ge \veps$. It follows that a subsequence $(y_n; n\in N)$ of $(x_n; n\in N)$
exists with  $d_i(y_n, x)\ge \veps$, $n\in N$; wherefrom, no convergent 
subsequence $(z_n; n\in N)$ of it (hence of $(x_n; n\in N)$) can have $x$ as limit, 
contradicting the above conclusion.	
\eproof

Concerning the notion we just introduced, it would be desirable 
(for both practical and theoretical reasons) to express it in terms of the sequence itself. 
To this end, let us call the sequence $(x_n; n\in N)$ in $X$, {\it precompact}
when for each $i\in N$, $\veps> 0$, a finite subset $A = A_{i,\veps}$  of $N$ may be found so that: 
for every $n\in N$ there exists $p\in A$ with $d_i(x_n, x_p)< \veps$. 
Now, as a completion of Lemma \ref{le1}, we have

\blemma \label{le2}
Assume $(X, D;\le)$ is such that $X$ is order complete. 
Then, for each increasing sequence in $X$, relatively compact is identical with precompact.
\elemma

\bproof
{\it Necessity}. Let $(x_n; n\in N)$ be an increasing relatively compact sequence
in $X$ that is not precompact. Then, a couple $i\in N$, $\veps> 0$ may be chosen so that,
for each finite subset $A$ of $N$, an index $n\in N$ will exist with $d_i(x_n,x_p)\ge \veps$, 
for all $p\in A$. It easily follows that a subsequence $(y_n; n\in N)$ of $(x_n; n\in N)$
may be constructed such that $d_i(y_n,y_m)\ge \veps$, $n<m$, proving $(y_n; n\in N)$ 
has no $D$-Cauchy (hence, by our hypothesis, no $D$-convergent) subsequences, 
contrary to our assumption.
{\it Sufficiency}. Let $(x_n; n\in N)$ be an increasing precompact sequence in $X$
and let $(y_n; n \in N)$ be a subsequence of it. As $(y_n; n \in N)$ is precompact too,
it clearly follows, by definition, that a subsequence $(u_n; n\in N)$ of it may be
found with $d_0(u_n, u_m)< 2^{-0}$, $n\le m$; furthermore, by the precompactness of
$(u_n; n\in N)$, a subsequence $(v_n; n\in N)$ of it may be found with 
$d_1(v_n, v_m)< 2^{-1}$, $n\le m$; and so on. 
By a standard diagonal process one arrives at an increasing $D$-Cauchy 
(hence, by our completeness hypothesis, $D$-convergent) 
subsequence $(z_n; n \in N)$ of $(y_n; n \in N)$ and the proof is complete.
\eproof

{\bf (B)}
Let $X$ be an ordered metrizable uniform space
under the denumerable sufficient family of semi-metrics $D=(d_i; i\in N)$ 
and the ordering $\le$. Also, let $Y$ be a subset of $X$ 
and $T$ a mapping from $Y$ to itself.
An important problem concerning these elements is that of determining
the existential comparative (modulo $\le$) connections between
the subset $Y_{oi}$ of all solutions in $Y$ of the operator inequality
\beq \label{401}
x\le Tx	
\eeq
and the subset $Y_{oe}$  of all solutions in $Y$ 
of the associated operator equation
\beq \label{402}
x=Tx.
\eeq
This will necessitate, as a first basic hypothesis
\bit
\item[(d01)]
$Y_{oi}$ is not empty.
\eit
In the following, we are interested in establishing 
a number of topological answers to the above question;
so, it is natural to accept as a second basic hypothesis
\bit
\item[(d02)]
$\le$ is interval-closed.
\eit
The first main result of the present paper is

\btheorem \label{t6}
Let the order-closed subset $Y$ of $X$ and 
the increasing mapping $T$ from $Y$ to itself be such that
\bit
\item[(d03)]
each increasing sequence $(x_n; n\in N)$ in $Y$ with $x_n\in T^{k(n)}(Y_{oi})$,
$n\in N$, for a strictly increasing sequence $(k(n); n\in N)$ in $N$, 
is relatively compact.
\eit
Then, to any $u$ in $Y_{oi}$ there corresponds 
$v$ in $Y_{oe}$ with the properties 

(a) $u\le v$,

(b) if $w\in Y_{oi}$ satisfies $v\le w$ then $v=w$.
\etheorem

\bproof
There are three steps to be passed.

{\bf I)}
Without loss of generality, 
one may suppose that $D$ is an increasing family ($d_i\le d_j$ whenever $i\le j$); 
because, otherwise, replacing it by the family 
$E=(e_i:=d_0+...+d_i; i\in N)$ 
the basic hypothesis (d02) as well as the specific assumption (d03) remain valid. 

{\bf II)}
We claim that for every couple $i\in N$, $\veps> 0$, 
the following assertion is true
\beq  \label{403}
\mbox{  \btab{l}
for each $m\in N$ and $x\in T^m(Y_{oi})$, there exist \\
$n> m$ in $N$ and $y\ge x$ in $T^n(Y_{oi})$  such that: \\
$d_i(y,z)< \veps$, 
for every $p> n$ in $N$ and $z\ge y$ in $T^p(Y_{oi})$.
\etab
}
\eeq
Indeed, if (\ref{403}) were not valid, 
a $m\in N$ and $x\in T^m(Y_{oi})$  exist with 
$$ 
\mbox{  \btab{l}
for every $n> m$ in $N$ and $y\ge x$ in $T^n(Y_{oi})$, there exist \\
$p> n$ in $N$ and $z\ge y$ in $T^p(Y_{oi})$ 
with $d_i(y,z)\ge \veps$.   
\etab
}
$$
It immediately follows that an increasing sequence $(y_n; n\in N)$ in $Y$ and a strictly
increasing sequence $(k(n); n\in N)$ in $N$ may be constructed with
$$
\mbox{
$y_n\in T^{k(n)}(Y_{oi})$\  and\   $d_i(y_n,y_{n+1})\ge \veps$,\  for all $n\in N$.
}
$$
By (d03), $(y_n; n\in N)$ is relatively compact;
hence $D$-convergent if we take (d02) plus Lemma \ref{le1} into account; 
so that $d_i(y_n, y_{n+1})\to 0$ as $n \to \oo$. 
The contradiction at which we arrived shows that the assertion (\ref{403}) is true.

{\bf III)}
For the arbitrary fixed $u$ in $Y_{oi}=T^0(Y_{oi})$ there exists, by (\ref{403}), 
an increasing sequence $(x_n; n\in N)$ in $Y$ and a strictly increasing sequence
$(k(n); n\in N)$ in $N$, fulfilling 
[$0< k(n)$, $u\le x_n\in T^{k(n)}(Y_{oi})$, $n\in N$], as well as
\beq \label{404}
\mbox{
($\forall n$):\ [$N\ni p> k(n)$  and  $T^p(Y_{oi})\ni y\ge x_n$]\ imply $d_n(y,x_n) < 2^{-n}$.
}
\eeq
From (d03), in conjunction with (d02) and Lemma \ref{le1}, 
it follows that $x_n \Dtends v$ for some $v$ in $Y$.
We claim that $v$ is the desired element. 
Indeed, let us first observe that, 
in view of the self-closeness property of our ordering,
\beq \label{405}
u\le x_n\le v,\ \ n\in N;	
\eeq
and therefore, $u\le v$. As an immediate consequence of (\ref{405}) 
we have $Tx_n \le Tv$, $n\in N$.
On the other hand, by the evident relation
$$
x_n \le Tx_n\in T^{k(n)+1}(Y_{oi}),\ \  n\in N
$$
plus (\ref{404}) it clearly follows $Tx_n\Dtends v$;
so, combining these, one arrives 
(by the anti-self-closeness property of our ordering) at $v\le Tv$; that is, $v\in Y_{oi}$. 
Moreover, as a further consequence of (\ref{405})
$$
x_n\le T^{k(n)+1}x_n\le T^{k(n)+1}v\in T^{k(n)+1}(Y_{oi}),\  \ n\in N;
$$
in which situation, again by (\ref{404}), $T^{k(n)+1}v\Dtends v$;
so that (via (d02))
$$
v\le Tv\le T^{k(n)+1}v\le v,\ \  n\in N;
$$
that is, $v\in Y_{oe}$. Finally, suppose that $v\le w$ for some $w\in Y_{oi}$;
then, as
$$
x_n\le v\le T^{k(n)+1}w \in T^{k(n)+1}(Y_{oi}),\ \  n\in N
$$
one gets by (\ref{404}) that $T^{k(n)+1}w \Dtends v$; and therefore, by (d02) again,
$$
\mbox{
$w\le T^{k(n)+1}w\le v,$  $n\in N$ [hence, in particular, $w\le v$],
}
$$
completing the argument.	
\eproof

Let us call the subset $Z$ of $X$, {\it order-sequentially relatively compact} 
when each increasing sequence 
in $Z$ is relatively compact.
Clearly, a sufficient condition guaranteeing the validity of (d03) is
\bit
\item[(d04)]
$T^k(Y)$ is order-sequentially relatively compact, for some index $k\in N$.
\eit
Then, as an useful variant of the first main result, we have 
(cf.
Turinici \cite{turinici-1984})

\btheorem \label{t7}
Let the order-closed subset $Y$ of $X$ and the increasing mapping $T$ 
from $Y$ to itself be such that (d04) is holding.
Then, conclusions (a)+(b) of Theorem \ref{t6} are retainable.
\etheorem

Returning to the setting of (d03) - essential to the present discussion - 
let us remark that its particular form (d04) 
may be viewed as a "spatial" (strong) restriction of it;
so that it is of practical interest to determine what happens when (d03) 
is replaced by its "temporal" (weak) restriction
\bit
\item[(d05)]
each increasing sequence $(T^nx; n\in N)$ in $Y$ with $x\in Y_{oi}$, \\
is relatively compact.
\eit
To do this, we have to introduce the notions below. 
Given the mapping $U$ from $Y$ to itself, call it 
{\it sequentially continuous at the left} when for each $x$ in $Y$ and 
each increasing sequence $(x_n; n\in N)$ in $Y$ with 
$x_n \Dtends x$ and $x_n\le x$, $n\in N$, we have $Ux_n\Dtends Ux$. 
Also, let us say that $U$ has an {\it order uniqueness property} 
when $x\le y$ and $x= Ux$, $y=Uy$ imply $x=y$ 
(i.e.: any two fixed points of $U$  are either identical or incomparable). 
The second main result of the present paper is 
(cf. also 
Dugundji and Granas \cite[Ch I, Sect 4]{dugundji-granas-1982})

\btheorem \label{t8}
Let the order-closed subset $Y$ of $X$ 
and the increasing mapping $T$ from $Y$ to itself be such that 
(d05) holds, as well as
\bit
\item[(d06)]
$T$ is sequentially continuous at the left
\item[(d07)]
$T$ has an order uniqueness property.
\eit
Then, conclusions (a)+(b) of the main result remain valid.
\etheorem

\bproof
Let $u$ in $Y_{oi}$ be arbitrary fixed. By (d05) plus (d02) and Lemma \ref{le1},
$T^nu\Dtends v$ for some $v\in Y$. Clearly, $T^nu\le v$, $n\in N$; so that, 
by the sequential left continuity assumption (d06), $T^{n+1}u\Dtends Tv$; 
wherefrom (as $D$ is sufficient), $v\in Y_{oe}$. 
Let $w$ in $Y_{oi}$ be such that $v\le w$. By the above reasoning $T^nw \Dtends v'$ 
for some $v'\in Y_{oe}$; on the other hand [by (d02)], $v\le T^nw\le v'$, $n\in N$, 
and this proves $v\le v'$. Combining this fact with (d07), one gets $v=v'$ 
and hence $w\le v$.
\eproof

An interesting feature of the above statements is given by the fact that 
(although implicitly embodied into the hypothesis (d03) or its variants)
no explicit (order) completeness property  for the ambient 
ordered metrizable uniform space were assumed; so that --
to complete our treatment and, at the same time, 
to cover some useful particular cases --
it would be necessary to discuss this eventuality.
Assume in the following that 
[in addition to the basic hypotheses (d01)+(d02)]
\bit
\item[(d08)]
$X$ is order complete;
\eit
then, in view of Lemma \ref{le2}, an appropriate
formulation of the main results might be obtained if one replaces 
in (d03), (d04), (d05), the word "relatively compact" by "precompact". 
Particularly, if we restrict our considerations to Theorem \ref{t8} above,
the following remark turns out to be in effect in many concrete situations.
Let $f:R_+\to R_+$ be increasing; 
we shall say that it has the property (P), provided
\bit
\item[]
$f^n(t)\to 0$\ \  as\ \  $n\to \oo$,\ \ for all $t>0$,
\eit
where, for each $n$, $f^n$ indicates the $n$-th iterate of $f$; 
note that, by a lemma due to 
Matkowski \cite{matkowski-1977},
we necessarily have in such a case $f(t)<t$, for all $t>0$ (hence $f(0)=0$).
Now, $Y$ and $T$ being as before, let us denote
\bit
\item[]
$f_i(t)=\sup\{d_i(Tx,Ty); x\le y,  d_i(x,y)\le t\}$,\ \ $t\in R_+$,\ $i\in N$.
\eit
Then we claim that the hypothesis
\bit
\item[(d09)]
$f_i(R_+)\incl R_+$ and $f_i$ has the property (P), for all $i\in N$  
\eit
is a sufficient one for the validity of (d05)+(d06)+(d07). 
Indeed, letting $u\in Y_{oi}$ be arbitrary fixed, put 
$a_i=d_i(u,Tu)$, $i\in N$, and observe that
$$
d_i(T^nu,T^{n+1}u)\le f_i^n(a_i),\ \  i,n\in N; 
$$
a relation which in turn implies, by (d09)
$$
\mbox{
$d_i(T^nu,T^{n+1}u)\to 0$\ \ as\ \  $n\to \oo$,\ \ for all $i\in N$.
}
$$
Let $i\in N$ and $\veps> 0$ be arbitrary fixed. 
By the above relation, a rank  $m=m(i,\veps)$ may be found with 
$d_i(T^mu,T^{m+1}u)\le \veps-f_i(\veps)\le \veps$;
combining with the definition of $f_i$ yields 
$d_i(T^{m+1}u,T^{m+2}u)\le f_i(\veps)$ 
so that, by the triangle property, $d_i(T^mu,T^{m+2}u)\le \veps$. 
Again using the definition of $f_i$ we have 
$d_i(T^{m+1}u,T^{m+3}u)\le f_i(\veps)$;
so that, by the same procedure as above, 
$d_i(T^mu,T^{m+3}u)\le \veps$, and so on.
By a finite induction, one arrives at 
$d_i(T^mu,T^{m+n}u)\le \veps$, $n\in N$.
This, along with (d08), proves (d05); so, the assertion follows,
because (d06)+(d07) are almost trivial in our case.
\sk  \sk

In concluding this section, let us remark that 
the comparison theorems we formulated before may 
be interpreted in the following dual ways:

{\bf i)}
as maximality results modulo $Y_{oi}$;
in which case, via Theorem \ref{t1} of 
Turinici \cite{turinici-1984} 
they appear  as particular versions of the maximality principle in 
\cite{turinici-1982}

{\bf ii)} 
as fixed point results modulo $Y$;
in which situation (under a continuity assumption similar to (d06))
they may be viewed as metrizable uniform versions of 
some topological statements in this area due to 
Wallace \cite{wallace-1945}, 
Ward \cite{ward-1954}, 
Smithson \cite{smithson-1971},
and 
Turinici \cite{turinici-1979}. 

\n
On the other hand, suppose that $X$ is a complete Fr\'{e}chet space 
under a denumerable sufficient family of seminorms $S=\{|.|_i; i\in N\}$ 
and let $X_+$ be a closed cone in $X$; 
then, defining an ordering structure by
\bit
\item[]
$x\le y$\ \  if and only if\ \  $y-x \in X_+$
\eit
the general hypotheses (d02)+(d08) of this section are clearly fulfilled; 
in particular, when $S$ reduces to a single element (i.e., a norm on X)
Theorem \ref{t6} includes the 1973 related statement in 
Krasnoselskii and Sobolev \cite{krasnoselskii-sobolev-1973}; 
and Theorem \ref{t8} reduces to the 1970 result in
Chandra and Fleishman \cite{chandra-fleishman-1970};
see also 
Azbelev and Tsaljuk \cite{azbelev-tsaljuk-1958}.
Some concrete examples of such cones may be found in 
Krasnoselskii \cite[Ch I]{krasnoselskii-1962}; 
cf. also 
Vulikh \cite[Ch III]{vulikh-1961}.
Finally, suppose the self-mapping $T$ were decreasing; 
then, evidently, $T^2$ is increasing; so that -- 
modulo the remaining hypotheses --
a number of comparison results concerning the couple (\ref{401})+(\ref{402})
(with $T^2$ in place of $T$) may be given. 
Some topological aspects of the problem were discussed by 
Seda \cite{seda-1981}.
\bk

{\bf (C)}
{\it Note added in 2011}
\sk

The argument concerning the couple (d08)+(d09) 
tells us that the following particular version 
of Theorem \ref{t8} was established.
(As before, (d01)+(d02) prevail).

\btheorem \label{t9}
Let the order-closed subset $Y$ of $X$ 
and the increasing mapping $T$ from $Y$ to itself be such that 
(d08)+(d09) hold.
Then, conclusions (a)+(b) of the main result are retainable.
\etheorem

In fact, a close examination of the reasoning above 
tells us that such conclusions are obtainable even if 
(d02) is to be replaced by its weaker counterpart 
\bit
\item[(d10)]
$\le$ is self-closed.
\eit
In this case, Theorem \ref{t9}
extends the 2008 statement in 
Agarwal, El-Gebeily and O'Regan \cite[Theorem 2.1]{agarwal-el-gebeily-o-regan-2008}
to the realm of ordered metrizable uniform spaces.



\begin{thebibliography}{99}


\bibitem{agarwal-el-gebeily-o-regan-2008}
{R. P. Agarwal}, {M. A. El-Gebeily} and {D. O'Regan},
\it Generalized contractions in partially ordered metric spaces, 
\rm Appl. Anal., 87 (2008), 109-116. 


\bibitem{altun-erduran-2011}
{I. Altun} and {A. Erduran},
\it A Suzuki type fixed-point theorem,
\rm Intl. J. Mathematics Math. Sci.,
Volume 2011, Article ID 736063, 9 pages, 2011.


\bibitem{azbelev-tsaljuk-1958}
{N. V. Azbelev} and {Z. B. Tsaljuk}, 
\it On the Chaplygin's problem
\rm [Russian], Ukrain. Mat. Zh., 10 (1958), 3-12. 


\bibitem{banach-1922}
{S. Banach}, 
\it Sur les op\'{e}rations dans les ensembles abstraits 
et leur application aux \'{e}quations int\'{e}grales.
\rm Fund. Math., 3 (1922), pp. 133-181.


\bibitem{chandra-fleishman-1970}
{J. Chandra} and {B. A. Fleishman}, 
\it On a generalization of the Gronwall-Bellman lemma in partially ordered Banach spaces,
\rm J. Math. Anal. Appl., 31 (1970), 668-681.


\bibitem{ciric-mihet-saadati-2009}
{L. B. Ciric}, {D. Mihet} and {R. Saadati},
\it Monotone generalized contractions in partially ordered probabilistic metric spaces,
\rm Topology and its Appl., 156 (2009), 2838-2844.


\bibitem{collaco-e-silva-1997}
{P. Collaco} and {J. C. E Silva},
\it A complete comparison of 25 contraction conditions,
\rm  Nonlin. Anal., 30 (1997), 471-476.


\bibitem{dugundji-granas-1982}
{J. Dugundji} and {A. Granas},
\it Fixed Point Theory,
\rm Vol. I, Monografie Mat., Vol 61, P.W.N., Warszawa, 1982.


\bibitem{gwozdz-lukawska-jachymski-2009}
{G. Gwozdz-Lukawska} and {J. Jachymski},
\it IFS on a metric space with a graph structure and extensions of the Kelisky-Rivlin theorem, 
\rm J. Math. Anal. Appl., 356 (2009), 453-463. 


\bibitem{kannan-1968}
{R. Kannan},
\it Some results on fixed points, 
\rm Bull. Calcutta Math. Soc., 60 (1968), 71-76. 


\bibitem{kikkawa-suzuki-2008}
{M. Kikkawa} and {T. Suzuki}, 
\it Some similarity between contractions and Kannan mappings, 
\rm Fixed Point Th. Appl., Volume 2008, Article ID 649749, 8 pages, 2008.


\bibitem{krasnoselskii-1962}
{M. A. Krasnoselskii}, 
\it Positive Solutions of Operator Equations
\rm [Russian], Gos. Izd. Fiz.-Mat. Lit., Moskva, 1962.


\bibitem{krasnoselskii-sobolev-1973}
{M. A. Krasnoselskii} and {A. V. Sobolev}, 
\it On the fixed points of non-continuous operators 
\rm [Russian], Sibirskii Mat. Zh., 14 (1973), 674-677.


\bibitem{matkowski-1977}
{J. Matkowski},
\it Fixed point theorems for mappings with a contractive iterate at a point,
\rm Proc. Amer. Math. Soc., 62 (1977), 344-348.


\bibitem{nieto-rodriguez-lopez-2007}
{J. J. Nieto} and  {R. Rodriguez-Lopez}, 
\it Existence and uniqueness of fixed point in partially ordered sets 
and applications to ordinary differential equations,
\rm Acta Math. Sinica (English Series), 23 (2007), 2205-2212.


\bibitem{o-regan-petrusel-2008}
{D. O'Regan} and  {A. Petru\c{s}el}, 
\it Fixed point theorems for generalized contractions in ordered metric spaces, 
\rm J. Math. Anal. Appl., 341 (2008), 1241-1252.


\bibitem{popescu-2009}
{O. Popescu}, 
\it Two fixed point theorems for generalized contractions with constants 
in complete metric space,
\rm Central European J. Math., 7 (2009), 529-538.


\bibitem{ran-reurings-2004}
{A. C. M. Ran} and  {M. C. Reurings}, 
\it A fixed point theorem in partially ordered sets 
and some applications to matrix equations, 
\rm Proc. Amer. Math. Soc., 132 (2004), 1435-1443.


\bibitem{rhoades-1977}
{B. E. Rhoades},
\it A comparison of various definitions of contractive mappings,
\rm Trans. Amer. Math. Soc., 226 (1977), 257-290.


\bibitem{rus-2001}
{I. A. Rus}, 
\it Generalized Contractions and Applications, 
\rm Cluj University Press, Cluj-Napoca, 2001.


\bibitem{seda-1981}
{V. Seda},
\it Antitone operators and ordinary differential equations. 
\rm Czech. Math. J., 31 (106) (1981), 531-553.


\bibitem{singh-pathak-mishra-2010}
{S. L. Singh}, {H. K. Pathak}, and {S. N. Mishra}
\it On a Suzuki type general fixed point theorem with applications,
\rm Fixed Point Th. Appl., Volume 2010, Article ID 234717, 15 pages, 2010.


\bibitem{smithson-1971}
{R. E. Smithson}, 
\it Fixed points of order preserving multifunctions. 
\rm Proc. Amer. Math. Soc., 28 (1971), 304-310.


\bibitem{suzuki-2008}
{T. Suzuki}, 
\it A generalized Banach contraction principle that characterizes metric completeness,
\rm Proc. Amer. Math. Soc., 136 (2008), 1861-1869.


\bibitem{turinici-1979}
{M. Turinici}, 
\it Abstract monotone mappings and applications to functional differential equations, 
\rm Atti Accad. Naz. Lincei (8), 66 (1979), 189-193.


\bibitem{turinici-1982}
{M. Turinici}, 
\it Constant and variable drop theorems on metrizable locally convex spaces, 
\rm Comment. Math. Univ. Carolin., 23 (1982), 383-398.


\bibitem{turinici-1984}
{M. Turinici}, 
\it Abstract Gronwall-Bellman inequalities on ordered metrizable uniform spaces, 
\rm J. Integral Equations, 6 (1984), 105-117.


\bibitem{turinici-1986}
{M. Turinici},
\it Abstract comparison principles and 
multivariable Gronwall-Bellman inequalities,
\rm J. Math. Anal.  Appl., 117 (1986), 100-127.


\bibitem{vulikh-1961}
{B. Z. Vulikh},
\it An Introduction to the Theory of Partially Ordered Spaces 
\rm [Russian], Gos. Izd. Fiz.-Mat. Lit., Moskva, 1961. 


\bibitem{wallace-1945}
{A. D. Wallace},
\it A fixed point theorem, 
\rm Bull. Amer. Math. Soc., 51 (1945), 413-416.


\bibitem{ward-1954}
{L. E. Ward, Jr.}, 
\it Partially ordered topological spaces,
\rm Proc. Amer. Math. Soc, 5 (1954), 144-161.



\end{thebibliography}
\end{document}